\documentclass[12pt]{article}

\usepackage{amssymb,latexsym,amsmath}
\usepackage{verbatim}
\newtheorem{theorem}{Theorem}
\newtheorem{lemma}[theorem]{Lemma}
\newtheorem{proposition}[theorem]{Proposition}

\newtheorem{corollary}[theorem]{Corollary}
\newtheorem{remark}[theorem]{Remark}

\newenvironment{proof}{\begin{trivlist}\item {\it
       Proof.\,}}{\mbox{}~\hfill~$\Box$\end{trivlist}}
\newcommand{\A}{\mathcal{A}}
\newcommand{\B}{\mathcal{B}}

\def\<{\langle}
\def\>{\rangle}

\begin{document}

\title{Normal Toric Ideals of Low Codimension}
\author{Pierre Dueck, Serkan Ho\c{s}ten and Bernd Sturmfels}
\date{}

\maketitle
\abstract{
Every normal toric ideal of codimension two is
minimally generated by a Gr\"obner basis with squarefree initial monomials.
A polynomial time algorithm is presented for checking  whether a
 toric ideal of fixed codimension is normal.
}

\section{Introduction}

Let $\A$ be a nonnegative integer $d \times n$ matrix of rank $d$,  
$C(\A)$ the cone in $\mathbb{R}^d$ spanned by the columns of $\A$, and  
 $I_{\A}$ the toric ideal  associated to $\A$
 as in \cite{gb-convex-polytopes}.
 The codimension of $I_\A$ is $n-d$.
If the columns of $\A$ form a Hilbert basis, i.e.~$\,C(\A) \cap \mathbb{Z} \A  = \mathbb{N} \A$, 
 then $\A$ and $I_{\A}$ are
 called {\em normal}.  Our main result is

\begin{theorem} \label{main-thm}
Every normal toric ideal $I_\A$ of codimension $2$ has a squarefree initial ideal,
and the corresponding reduced Gr\"obner basis minimally generates~$I_\A$.
\end{theorem}

In Section 2 we introduce the necessary tools and prove Theorem \ref{main-thm} in the case of complete intersections.  In Section 3 we complete the proof in the case of codimension two 
normal toric ideals which are not complete intersections. Section 4 deals with 
detecting normality and we prove the following result.

\begin{theorem} \label{complexity}
If $ \mathrm{codim}(I_\A) = n{-}d$ is fixed, then there is an
algorithm to decide whether $\A$ is normal whose running time is polynomial
in $n$ and the bitsize of~$\A$.
 \end{theorem}    

In other words, if $n-d$ is fixed then one can decide in polynomial time whether $n$ vectors in 
$\mathbb{Z}^d$ form a Hilbert basis of the cone these vectors generate.
This answers a question raised by Alexander Barvinok at 
Snowbird (June 2006). We recently found out that Theorem \ref{complexity} has been proved  independently by F.~Eisenbrand, A.~Seb\"o, and G.~Shmonin \cite{ESS}. The authors 
had announced the result at the 12th Combinatorial Optimization Workshop held in Aussois, France 
in January 2008. For recent advances on related
complexity questions concerning lattice points in polyhedra
of fixed dimension or codimension we refer to the article \cite{ES}.

\smallskip

This work is motivated both by questions that are intrinsic to 
combinatorial commutative algebra and by 
applications to statistics and optimization.
In the former domain,
a longstanding conjecture states that every Cohen-Macaulay toric ideal $I_{\A}$ 
has a monomial initial ideal that is also Cohen-Macaulay. 
This holds for toric ideals up to
dimension 3 (see \cite{gomory}), and 
O'Shea and Thomas proved it for $\Delta$-normal configurations   \cite{delta-normal}.
Theorem \ref{main-thm} offers supporting evidence because normal toric ideals 
and their square-free initial ideals are both Cohen-Macaulay.  

In algebraic statistics, our result ensures that the
sequential importance sampling  scheme of 
Chen, Dinwoodie and Sullivant \cite{CDS}
is applicable to exponential families with few states.
In integer programming, it ensures that,
for suitably chosen cost functions,
  every matrix of corank $2$ specifies a Gomory family \cite{gomory}.
Finally, in algebraic geometry, where 
the definition of toric varieties  \cite{Fulton} requires them to be normal,
Theorem \ref{main-thm} states that
every toric variety of codimension 2 admits
a Gr\"obner degeneration to a reduced union
of coordinate subspaces.
 
\section{Complete Intersections}

The aim of this section is to prove Theorem  \ref{main-thm}
in the case when $I_\A$ is a complete intersection.
Before we approach this proof we will prove an alternative
characterization of normality due to Seb\"o \cite{sebo}.
We  include a proof of this result.

\begin{proposition} \label{Hilbert-kernel} An integer matrix 
$\A = [a_1,\ldots,a_n]$ 
is normal if and only if for each $x \in \ker(\A)$ there exists an integer vector $y \in \ker(\A)$ such that $y \leq \lceil x \rceil$.
\end{proposition}

\begin{proof} ($\Longrightarrow$)
Suppose $\mathcal{A} $ is normal, and $x \in \ker \A$.  We have
$\,x_1a_1 + \cdots + x_n a_n = 0$.
The vector $\,z \, = \, \lceil x_1 \rceil \cdot a_1 +  \,\cdots\, + \lceil x_n \rceil  \cdot a_n \,$
lies in the lattice $\mathbb{Z} \A$, and as
$\,z \, = \,( \lceil x_1 \rceil - x_1 ) \cdot a_1 + \,\cdots\, + 
(\lceil x_n \rceil -x_n) \cdot a_n $, it also lies in the cone $C(\A)$.
Since $\A$ is normal, we conclude that $z$ is in
the semigroup $\mathbb{N} \A$. We can write
$ z = m_1 a_1 + \cdots + m_n a_n $ with $m_1,\ldots,m_n$
nonnegative integers. The vector $y $
with coordinates  $y_i = \lceil x_i \rceil - m_i$ lies in 
$ \ker( \A)$ and satisfies $y \leq \lceil x \rceil $.   \smallskip

($\Longleftarrow$) Now suppose that for each $x \in \ker (\A)$ there is an 
integral $y \in \ker (\A)$ with $y \leq \lceil x \rceil$.  Let $z \in C(\A) \cap \mathbb{Z}^d$.  This means that
$  z = r_1 a_1 + \ldots + r_n a_n $ with $ r_i \in \mathbb{R}_{\geq 0}  $ and 
$  z = m_1 a_1 + \ldots + m_n a_n $ with $ m_i \in \mathbb{Z} . $
Combining these we obtain:
\[ (m_1 - r_1) \cdot a_1 \,+ \, \cdots \,+\, (m_n - r_n) \cdot a_n = 0.\]
By hypothesis we may pick an integral $y \in \ker (\A)$ with 
$ y_i \leq \lceil m_i - r_i \rceil $ Then
\[ z \,\,\,= \,\,\, (m_1 - y_1) \cdot a_1 \,+\, \cdots \,+\, (m_n - y_n) \cdot a_n \]
gives a nonnegative integral representation of $z$ in terms of columns of $\A$,
since $m_i -y_i \geq 0$ for all $i$.
We conclude that $z \in \mathbb{N} \A$, and hence $\A$ is normal. 
\end{proof}

Our main tool in what follows
is the Gale diagram of a vector configuration.  Let $\A$ be a integer $d \times n$ matrix whose column vectors span $\mathbb{R}^d$.  We choose a matrix $\B$ whose rows form a lattice basis of $\ker(\A) \cap \mathbb{Z}^n$. The set of column vectors of $\B$ is said to be a {\em Gale diagram} \cite{polytopes-ziegler} of $\A$.  Normality of $I_{\A}$ is encoded in both $\A$ and in $\ker(\A) = {\rm im}(\B)$, by Proposition \ref{Hilbert-kernel}, and hence also in $\B$.

Hochster \cite{hochster} proved that a normal toric ideal is Cohen-Macaulay.  Thus a codimension two normal toric ideal has a minimal free resolution of length two.  
Peeva and Sturmfels  \cite{codim-2-lattice}  characterized Cohen-Macaulay codimension two lattice ideals.  In this paper we consider saturated lattices whose lattice ideal is a toric ideal.  For the remainder of this paper we assume that $I_{\A}$ is normal, or equivalently that $\A$ is a Hilbert basis of the cone $C(\A)$. We assume that the cone $C(\A)$ is pointed and hence $I_{\A}$ is 
homogeneous in some positive grading.
 
\smallskip

The following result gives a supply of squarefree monomial terms of the binomial generators of $I_{\A}$.   
This result has been proven in \cite[Proposition 4.1]{SV} and \cite[Lemma 6.1]{OH}, and we have 
also learned it from Winfried  Bruns \cite{B}.

\begin{proposition} \label{normal-minimal}
Suppose $\A$ is normal.  Then each minimal binomial generator of the toric ideal
 $I_{\A}$ has at least one squarefree term.  
\end{proposition}

This implies that the conclusion of Theorem 
\ref{main-thm} holds when ${\rm codim}(I_\A)= 1$. In that case, 
$I_A$ is a principal ideal and the unique binomial generator of $I_\A$ is a Gr\"obner basis
with its squarefree term being the leading monomial.

In view of Proposition~\ref{normal-minimal} our approach is to show the existence of a term order selecting the squarefree terms as initial terms. The Gale diagram gives information toward this goal.   
The following result is \cite[Proposition 4.1]{codim-2-lattice}.

\begin{proposition}
If $\,{\rm codim}(I_\A)= 2\,$ then the following are equivalent: \\
(i) The toric ideal $I_\A$ is not Cohen-Macaulay. \\
(ii) The toric ideal $I_\A$ has at least four minimal generators.  \\
(iii) The matrix $\A$ has a Gale diagram $\B$ which intersects each of the four open quadrants
in $\mathbb{R}^{2}$. Here the matrix $\B$ is identified with its set of columns.
\end{proposition}

\begin{figure}
\centering
\setlength{\unitlength}{2 mm} 
\begin{picture}(60,40) 
\thicklines
\put(30,20){\vector(1,0){10}}\thicklines
\put(30,20){\vector(1,1){8}} \thicklines 
\put(30,20){\vector(0,1){10}} \thicklines 
\put(30,20){\vector(-1,1){8}} \thicklines 
\put(30,20){\vector(-1,0){10}} \thicklines 
\put(30,20){\vector(0,-1){10}} 
\put(15.0,19.0){$E$}
\put(43.0, 19.0){$A$}
\put(16.0, 30.0){$D$}
\put(29.0, 33.0){$C$}
\put(42.0, 30.0){$B$}
\put(29.0, 7.0){$F$}
\end{picture}
\vskip -1.6cm
\[
\begin{array}{ccccccc}
 & A & B & C & D & E & F \\
p-q \,\,\,=  & + & + & 0 & - & - & 0 \\
r-s \,\,\, =   & 0 & + & + & + & 0 & - \\
\end{array}
\] 
\vskip -0.5cm
\caption{The imbalanced Gale diagram of a complete intersection.}
\label{fig-compinter}
\end{figure}

We now assume that $I_\A$ is normal of codimension $2$
and $\B$ is any Gale diagram. Then $\,\B = \bigl\{ ( \B_{1j}, \B_{2j})\,:\,
j=1,\ldots,n \bigr\}\,$
intersects at most three open quadrants, and that any minimal generating set of $I_{\A}$ has two or three elements.  In this section we examine the first case,
where $\,I_{\A} = \< x^p - x^q, \,\, x^r - x^s \>$ is a complete intersection, and $\B$ is the $2 \times n$ matrix whose rows are $ p-q$ and $r-s$.  The Gale diagram $\B$ is  called
 {\em imbalanced} if either $\B_{1j} = 0$ or $\B_{2j} \geq 0$ for all $j$.

\begin{lemma} \label{lem-imbalanced} {\rm \cite[Theorem 3.9]{codim-2-lattice}}
\ \
A codimension $2$ toric ideal $I_A$ is a complete intersection if and only if 
there exists an imbalanced Gale diagram $\B$.
\end{lemma}

In the light of this lemma, we can represent a complete intersection
by an imbalanced Gale diagram as depicted in Figure \ref{fig-compinter}.
  The arrows represent a {\em sign class} of columns of $\B$ and not just an individual vector. For instance the class labeled $D$ in  Figure \ref{fig-compinter}
 consists of all column vectors with $\B_{1j} < 0$ and $\B_{2j} > 0$.  

\medskip

\noindent
{\sl Proof of Theorem \ref{main-thm} for complete intersections: }
    By Proposition \ref{normal-minimal} both generators $g_1 = x^p - x^q$ and $g_2 = x^r - x^s$ have a squarefree term. The class of vectors $F$ must exist in the Gale diagram since otherwise $C(\A)$ would not be pointed. 
        
    If the term $x^s$ corresponding to $F$ is squarefree then we can use any term order so that $x^s$ is the initial term of $g_2$ and
      the squarefree term of $g_1$ is its initial term.
Then $g_1$ and $g_2$ forms the Gr\"obner basis of $I_\A$ because their
initial terms are relatively prime. This Gr\"obner basis is or can be made reduced.

Suppose  that
 $x^s$ is not squarefree. Then $-\B_{2j} =  f \geq 2$ for some $j \in F$,
        and $\B_{2j} = 1$ for $j \in B \cup C \cup D$.
Without loss of generality we assume that $x^p$ is
           the squarefree term of $g_1$, so that
           $\B_{1j} = 1$ for $j \in A \cup B$.
            We choose representatives from the $D$ and $E$ classes, labeling them $-d$ and $-e$ where $d,e \geq 1$:
$$
\begin{array}{cccccccccccc}
            & A & & B& & C& & D&   & E&  & F \\
p-q \, \,\,= & 1&\ldots & 1&\ldots & 0&\ldots  & - d&\ldots & - e&\ldots & 0 \\
r-s  \,\,\, = & 0 &\ldots& 1&\ldots & 1&\ldots  & 1&\ldots   & 0&\ldots   & - f \\
\end{array}
$$
Now we consider $ u = - \frac{1}{2} (p-q) + \frac{1}{2} (r-s) \in \ker(\A)$ and we
round it up to get
$$
\begin{array}{cccccccccccc}
            & A && B& & C& & D&   & E& & F \\
\lceil u \rceil \,\, = &0  & \!\ldots \!& 0 & \!\ldots\! & 1& \! \ldots \! & \left\lceil 
\frac{d+1}{2} \right\rceil & \ldots& \left\lceil \frac{e}{2} \right\rceil 
& \ldots & \left\lceil - \frac{f}{2} \right\rceil  \\
\end{array}
$$
Note that $- f/2 \leq -1$.  By Proposition \ref{Hilbert-kernel} there exists an integral $v \in \ker(A)$ with $v \leq \lceil u \rceil$.  
This vector is an integral combination $v = \alpha (p-q) + \beta (r-s)$.

 If $v_1 =0$ then $\alpha = 0$ and $v$ must be a \emph{positive} multiple of $(r-s)$ to ensure that $v_F = -\beta f \leq \lceil -f/2 \rceil \leq -1  $.  This implies the contradiction $v_B = \beta > 0$.

Next suppose that $v_1 \leq -1$. Then  $\alpha \leq -1$ (considering the $A$ component) and $\beta \geq 1$ (considering the $F$ component).  The $D$ representative requires that 
$d+1 \leq - \alpha d  +\beta \leq \left\lceil \frac{d+1}{2} \right\rceil$.  But this implies that $d=0$, a contradiction. 
We conclude that the $D$ class is not present in the Gale diagram. By rotating the
diagram by $90$ degrees counterclockwise we can assume that 
we have an imbalanced Gale diagram where the $B$ class is missing. For the two new minimal
generators $x^p - x^q$ and $x^r - x^s$ we are either in the first case analyzed above 
(i.e. $x^s$ and $x^p$ are squarefree and relatively prime)
or in the second case where $x^p$ and $x^r$ are squarefree and relatively prime,
as these do not contain any $B$ variable.  In both cases the
two generators form a  squarefree Gr\"obner basis.
 \hfill $\Box$ 
 
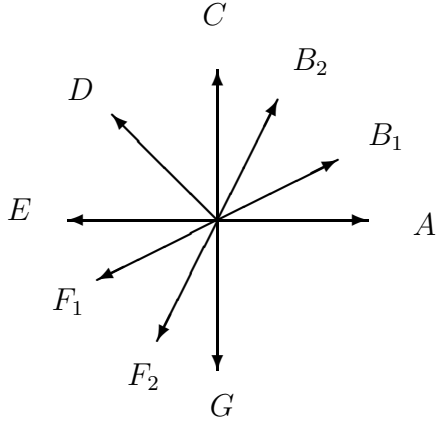
\begin{figure}[h]
\vskip -1cm
\centering
\setlength{\unitlength}{2 mm} 
\begin{picture}(60,40) 
\thicklines
\put(30,20){\vector(1,0){10}} 
\put(30,20){\vector(1,2){4}}
\put(30, 20){\vector(2,1){8}}
\put(30,20){\vector(0,1){10}} 
\put(30,20){\vector(-1,1){7}} 
\put(30,20){\vector(-1,0){10}} 
\put(30,20){\vector(0,-1){10}} 
\put(30,20){\vector(-1,-2){4}} 
\put(30, 20){\vector(-2, -1){8}}
\put(16.0,20.0){$E$} 
\put(43.0, 19.0){$A$} 
\put(20.0, 28.0){$D$} 
\put(29.0, 33.0){$C$}
\put(35.0, 30.0){$B_2$}
\put(40.0, 25.0){$B_1$}
\put(29.5, 7.0){$G$}
\put(19.0, 14.0){$F_1$}
\put(24.0, 9.0){$F_2$}
\end{picture}
\vskip -1.5cm
\caption{Gale diagram of a CM but not complete intersection configuration}
\label{fig-non-compinter}
\end{figure}

\begin{figure}
\vskip -1.8cm
\centering
\setlength{\unitlength}{2 mm}
\begin{picture}(40,50) 
\linethickness{0.075mm} 

\thicklines
\put(15,15){\circle*{1}}
\put(30, 15){\circle*{1}}
\put(15,30){\circle*{1}}
\put(15,15){\line(0,1){15}}
\put(15,15){\line(1,0){15}}
\put(15, 30){\line(1,-1){15}}

\put(15,30){\vector(1,2){3}}
\put(15,30){\vector(0,1){7}}
\put(15,30){\vector(-1,1){5}}

\put(15,15){\vector(-1,0){6}}
\put(15,15){\vector(-2,-1){6}}
\put(15,15){\vector(-1,-2){3}}
\put(15,15){\vector(0,-1){6}}

\put(30,15){\vector(1,0){6}}
\put(30,15){\vector(2,1){6}}

\put(38,15){$A$}
\put(38,18){$B_1$}
\put(19,37){$B_2$}
\put(14,39){$C$}
\put(7,36){$D$}
\put(6,15){$E$}
\put(6,10){$F_1$}
\put(10,6){$F_2$}
\put(14.5,6){$G$}
\put(16,16){$(0,0)$}
\put(16, 30){$(0,1)$}
\put(28, 12){$(1,0)$}
\end{picture}
\vskip -1cm
\caption{A syzygy triangle as in \cite{codim-2-lattice} }
\label{figure:first-minimal-syzygy}
\end{figure}
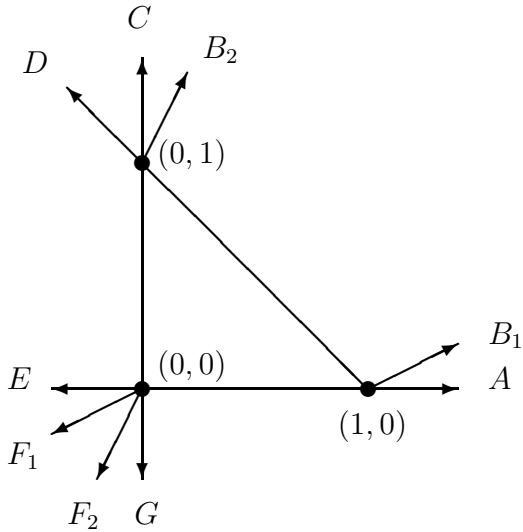

\section{Normal but not complete intersection}

We now assume that the toric ideal $I_\A$ is not a complete intersection, but it is 
normal and hence Cohen-Macaulay. This time we can assume that the Gale diagram is of the form 
as in Figure \ref{fig-non-compinter}. As in Section~2, the vectors in the
diagram represent a sign class of vectors. 
The minimal free resolution of $I_A$ has the form
\[ 0 \to R^2 \to R^3 \to R \to R/I_{\A} \to 0 \]  
where $R = \mathbb{K}[x_1, \ldots, x_n ] $. 
A matrix representing the map $R^2 \rightarrow R^3$
in the resolution can be determined using the
recipe in \cite{codim-2-lattice}. More precisely, using the
 two {\em syzygy triangles} in Figure 3 and Figure 4
we find that this matrix equals
\[
\begin{bmatrix}
AB_1 & D^*E^*F_1^* \\
B_2CD & F_2^*G^* \\
F_1F_2 & B_1^*B_2^* \\
\end{bmatrix}
\]

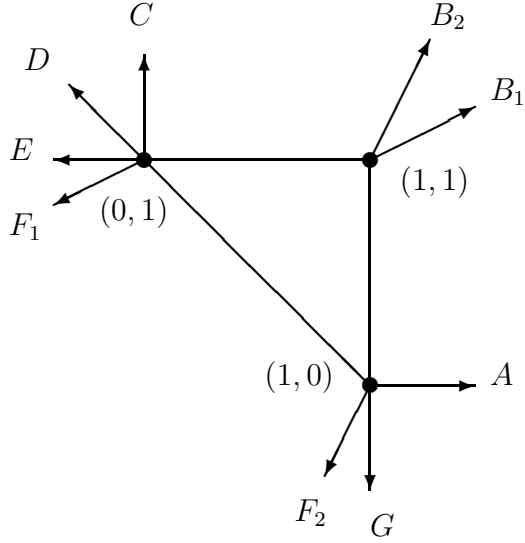
\begin{figure}
\vskip -.5cm
\centering
\setlength{\unitlength}{2 mm}
\begin{picture}(40,50) 
\linethickness{0.075mm} 

\thicklines
\put(30,30){\circle*{1}}
\put(30, 15){\circle*{1}}
\put(15,30){\circle*{1}}
\put(15,30){\line(1,0){15}}
\put(30,15){\line(0,1){15}}
\put(15, 30){\line(1,-1){15}}

\put(30,30){\vector(2,1){7}}
\put(30,30){\vector(1,2){4}}

\put(15,30){\vector(-1,1){5}}
\put(15,30){\vector(0,1){7}}
\put(15,30){\vector(-2,-1){6}}
\put(15,30){\vector(-1,0){6}}

\put(30,15){\vector(1,0){7}}
\put(30,15){\vector(-1,-2){3}}
\put(30,15){\vector(0,-1){7}}

\put(38,15){$A$}
\put(38,34){$B_1$}
\put(34,39){$B_2$}
\put(14,39){$C$}
\put(7,36){$D$}
\put(6,30){$E$}
\put(6,25){$F_1$}
\put(25,6){$F_2$}
\put(30,5){$G$}
\put(32,28){$(1,1)$}
\put(12, 26){$(0,1)$}
\put(23, 15){$(1,0)$}
\end{picture}
\vskip -1cm
\caption{Another syzygy triangle as in \cite{codim-2-lattice} }
\label{figure:second-minimal-syzygy}
\end{figure}

In the above matrix the letters represent classes of variables corresponding
to the classes of vectors in the Gale diagram. Products of letters correspond 
to monomials in these classes of variables. The letters with an asterisk
correspond to the same class as those without an asterisk, but
they might have different exponent vectors (since they come from
different syzygy triangles). The first column of the matrix corresponds to
the first triangle and the second column to the second triangle. 
Moreover, by the Hilbert-Burch Theorem, the three $2 {\times} 2$-minors of this matrix 
are precisely the minimal generators of $I_{\A}$:

\begin{equation} \label{bin1}
AB_1F_2^*G^* - B_2CDD^*E^*F_1^* 
\end{equation}
\begin{equation} \label{bin2}
AB_1B_1^*B_2^* - D^*E^*F_1F_1^*F_2 
\end{equation}
\begin{equation} \label{bin3}
B_1^*B_2B_2^*CD - F_1F_2F_2^*G^*
\end{equation}
The next result completes the proof of Theorem \ref{main-thm}.

\begin{lemma} \label{non-compinter}
There exists a term order such that the
binomial generators {\rm (\ref{bin1}),(\ref{bin2}),(\ref{bin3})}
of the toric ideal $I_\A$ form a Gr\"obner basis with squarefree initial monomials.
\end{lemma}

\begin{proof}
We have a few cases to consider. First, either the $D$ class
exists or it does not.
If it exists then the monomial $D D^*$ is not squarefree
and hence the squarefree term of (\ref{bin1})
is the first term. Note that the first terms of 
(\ref{bin2}) and (\ref{bin3}) cannot be squarefree simultaneously: if they were,
the $B_1$ and $B_2$ vectors cannot be present in the Gale diagram, and
this would be an imbalanced Gale diagram. Similarly, the second terms of 
these binomials cannot be squarefree simultaneously.
This gives two cases to consider. In the first case we have
\begin{equation}
\label{GBunderlined1}
 \underline{AG^*} - B_2CDD^*E^*F_1^*,  \,\,  \underline{AB_2^*} - D^*E^*F_1F_1^*, \, \,
B_2B_2^*CD - \underline{F_1G^*}.
\end{equation}
Here $B_1$ and $F_2$ are absent because otherwise
$B_1 B_1^*$ and $F_2F_2^*$ are not squarefree.
If we choose a lexicographic term order where $A > G > \{B_2, C, D, E, F_1\}$
then the underlined terms are the leading terms in (\ref{GBunderlined1}).
The S-pair $S(1,2) =  D^*E^*F_1F_1^*G^* - B_2B_2^*CDD^*E^*F_1^*$ 
is reduced to zero by the third binomial, and the S-pair 
$S(1,3) = AB_2B_2^*CD - B_2CDD^*E^*F_1F_1^*$ is reduced to zero
by the second binomial. The S-pair $S(2,3)$ reduces to zero since
the leading terms $AB_2^*$ and $F_1G^*$ are
relatively prime, and hence (\ref{GBunderlined1}) is a  squarefree Gr\"obner basis.

In the second case, the minimal generators and their squarefree terms are
\begin{equation}
\label{GBunderlined2}
\underline{AB_1F_2^*G^*} - CDD^*E^*,\,\,
AB_1B_1^* - \underline{D^*E^* F_2} , \,\,
\underline{B_1^* CD} - F_2F_2^*G^*.
\end{equation}
The product of the three underlined terms is equal to the
product of the three non-underlined terms. Hence
no term order selects the underlined terms as leading terms.
However, the squarefreeness of these three monomials implies
$$ A = \left[ \begin{array}{c} 1 \\ 0 \end{array} \right], \,
B_1 = \left[ \begin{array}{c} 1 \\ 1 \end{array} \right], \,
C = \left[ \begin{array}{c} 0 \\ 1 \end{array} \right], \,
D = \left[ \begin{array}{c} -1 \\ 1 \end{array} \right], $$ $$
E = \left[ \begin{array}{c} -1 \\ 0 \end{array} \right], \,\,\,
F_2 = \left[ \begin{array}{c} -1 \\ -1 \end{array} \right], \,\,\,
G = \left[ \begin{array}{r} 0 \\ -1 \end{array} \right].$$
From the diagonal edge in the two syzygy triangles we see that $B_1 = F_2^* = 1$.
This means that the non-underlined terms of the second and third binomials
in (\ref{GBunderlined2})
are actually squarefree, and we are back in the previous case
(\ref{GBunderlined1}).

Now suppose that the $D$ vectors are not present in the Gale diagram
depicted in Figure \ref{fig-non-compinter}.
Then the binomial generators (\ref{bin1}), (\ref{bin2}) and (\ref{bin3})  have the form
$$ AB_1F_2^*G^* - B_2CE^*F_1^*, \,\, AB_1B_1^*B_2^* - E^*F_1F_1^*F_2, \,\,
B_1^*B_2B_2^*C - F_1F_2F_2^*G^*.$$
If in the first binomial the first term is squarefree we are back to 
(\ref{GBunderlined1}) or (\ref{GBunderlined2}).
If the second term is squarefree, then we rotate the Gale diagram 
$180$ degrees. This leads to the same binomials but now with 
the first term of the first binomial squarefree. Once again we are
back to (\ref{GBunderlined1}) or (\ref{GBunderlined2}).
This concludes the proof.
\end{proof}

\section{Checking Normality}
In this section we assume that the codimension $m= n-d$ of $I_\A$ is fixed. First we reformulate 
Proposition \ref{Hilbert-kernel}. Let $z$ be an integral vector in $\mathrm{ker}(\A)$. 
We define 
$$ P_z \, \, = \, \, \bigl\{ \,x \in \mathrm{ker}(\A) \, : \, \lceil x_i \rceil \geq z_i \,\,\, 
\hbox{for} \,\, i=1, \ldots, n  \big\}.$$
Since $ P_z \, \, = \, \, \{ x \in \mathrm{ker}(\A) \, : \,  x_i > z_i-1 \,\,\, i=1, \ldots, n \}$
and since we assume that the cone $C(\A)$ is pointed, $P_z$ is a
relatively open polytope in
$\mathrm{ker}(\A) \simeq \mathbb{R}^m$.

\begin{remark} \label{shift}
If $u$ and $z$ are lattice vectors in $\mathrm{ker}(\A)$ then $P_{u+z} = 
z + P_u$.
\end{remark}

Now let $\B$ be an $n\times m$ matrix whose 
columns form a lattice basis of $\mathrm{ker}(\A)$, and let $b_i$ be the rows of $\B$. Then $P_z$ is affinely  isomorphic  to $Q_v = \{y \in \mathbb{R}^m \, \, : \, \, b_i \cdot y > b_i \cdot v - 1, \,\, i=1,\ldots,n\}$
where $v$ is the unique lattice point in $\mathbb{Z}^m$ such that $\B v=z$. Remark
\ref{shift} implies that for two lattice points $v$ and $w$ in $\mathbb{Z}^m$ we have 
$Q_{v + w} = w + Q_v$. Note that $Q_0 = \{y \in \mathbb{R}^m \, \, : \, \, b_i \cdot y > -1, \,\, i= 1, \ldots, n\}$.
We now see that the following is equivalent to Proposition \ref{Hilbert-kernel}.

\begin{theorem} The toric ideal $I_\A$ is normal if and only if $Q_0 + \mathbb{Z}^m = \mathbb{R}^m$.
\end{theorem} 

Given any polytope $Q = \{y \in \mathbb{R}^m \,\, : \,\, Cy \geq d \}$ of dimension $m$ the smallest 
positive real number $t$ such that $tQ + \mathbb{Z}^m = \mathbb{R}^m$ is called the 
{\it covering radius} of $Q$. If $Q$ is a rational polytope it is known that 
the covering radius of $Q$ is a rational number with a bit-size that is a polynomial in the bit-size 
of $C$ and $d$. 

\begin{corollary} The toric ideal $I_\A$ is normal if and only if the 
covering radius of $\bar{Q}_0$, the closure of the polytope $Q_0$, is less than $1$.
\end{corollary} 

\vskip 0.2cm 
\noindent {\it Proof of Theorem \ref{complexity}}: Ravi Kannan \cite[Section 5]{K} has shown that, for fixed $m$, and
given a rational $m$-dimensional polytope $Q = \{y \in \mathbb{R}^m \,\, : \,\, Cy \geq d \}$ where 
$C \in \mathbb{Z}^{n \times m}$ and $d \in \mathbb{Z}^n$, there exists an algorithm to find the 
covering radius of $Q$ with runtime a polynomial in $n$ and the bit-size of $C$ and the vector $d$.
Since one can compute a $\B$ whose bit-size is a polynomial in the bit-size of $\A$ in polynomial time,
the above corollary implies the result.  \hfill $\Box$

\bigskip
\medskip

\noindent
{\bf Acknowledgement}: B.~Sturmfels was supported
by the NSF  (DMS-0456960).

\bigskip

\noindent
{\em Authors' e-mail addresses}: 
{\tt dueck@math.ucdavis.edu}, \hfill \break
{\tt serkan@math.sfsu.edu},
{\tt bernd@math.berkeley.edu}

\end{document}